\date{} 
\title{Hyperfinite graph limits}
\author{Oded Schramm}
\documentclass[12pt]{article}
\usepackage{amsmath}
\usepackage{amssymb}
\usepackage{amsthm}
\usepackage{amsfonts}
\usepackage{graphicx}

\newif\iffigures\figurestrue
\figuresfalse

\newif\ifhyper\IfFileExists{hyperref.sty}{\hypertrue}{\hyperfalse}

\ifhyper\usepackage{hyperref}
\def\hitem#1#2{\item[\hypertarget{#1}{#2}]\expandafter\gdef\csname LBL#1ITM\endcsname{#2}}
\def\iref#1{\hyperlink{#1}{\csname LBL#1ITM\endcsname}}
\else
\def\hitem#1#2{\item[{#2}]\expandafter\gdef\csname LBL#1ITM\endcsname{#2}}
\def\iref#1{{\csname LBL#1ITM\endcsname}}
\fi

\newif\ifdraft
\drafttrue
\long\def\note#1/{\ifdraft {\bf [#1]}\fi}
\long\def\comment#1{}
\long\def\old#1{}

\numberwithin{equation}{section}
\numberwithin{figure}{section}

\newtheorem{theorem}{Theorem}
\numberwithin{theorem}{section}

\newtheorem{lemma}[theorem]{Lemma}

\theoremstyle{remark}
\theoremstyle{remark}

\def\LiptonTarjanApplications{MR584516}
\def\BSplanarLimits{MR1873300}
\def\AldousLyons{2006math......3062A}
\def\LovaszSzegedy{MR2274085}
\def\ElekCost{2006math......8474E}
\def\ElekLtwo{2007arXiv0709.1261E}
\def\ElekRegularity{arXiv:0711.2800}
\def\BSSpropertyTesting{BSSpropertyTesting}
\def\BSS{\BSSpropertyTesting}
\def\BSpyond{MR97j:60179}

\let\qqed=\qed
\def\QED{\qqed\medskip}
\let\qed=\QED

\newcommand{\N}{\mathbb{N}}

\def\SLEkk#1/{$\mathrm{SLE}(#1)$}
\def\SLEr#1/{$\mathrm{SLE(\kappa;#1)}$}
\def\SLEkr#1;#2/{$\mathrm{SLE(#1;#2)}$}
\def\SLEk/{\SLEkk{\kappa}/}
\def\SLEtwo/{\SLEkk2/}
\def\SLE/{$\mathrm{SLE}$}
\def\SLEab/{\SLEkr 4; {a/\hco-1}, {b/\hco-1}/}

\def\Ito/{It\^o}
\def \eps {\epsilon}
\def \P {{\bf P}}
\def\md{\mid}
\def\Bb#1#2{{\def\md{\bigm| }#1\bigl[#2\bigr]}}
\def\BB#1#2{{\def\md{\Bigm| }#1\Bigl[#2\Bigr]}}

\def\Pb{\Bb\P}
\def\Eb{\Bb\E}

\def\EB{\BB\E}

\def \p {{\partial}}
\def \E {{\bf E}}

\def\ev#1{{\mathcal{#1}}}
\def \proof {{ \medbreak \noindent {\bf Proof.} }}
\def\proofof#1{{ \medbreak \noindent {\bf Proof of #1.} }}

\def\UU{\mathcal{U}}

\def\bl{\bigl}\def\br{\bigr}\def\Bl{\Bigl}\def\Br{\Bigr}

\def\G{\widehat{\mathcal G}}
\def\GG{\G_M}
\def\GfinM{\mathcal G^0_M}
\def\MM{\mathfrak M_M}
\def\M{\mathfrak M}

\def\mB{\mathcal B}

\def\noopsort#1{}

\begin{document}
\maketitle

\begin{abstract}
G\'abor Elek introduced the notion of a hyperfinite graph family:
a collection of graphs is hypefinite if for every $\eps>0$ there is some finite $k$
such that each graph $G$ in the collection can be broken into connected components of size at most $k$
by removing a set of edges of size at most $\eps\,|V(G)|$.
We presently extend this notion to a certain compactification of finite
bounded-degree graphs, and show that
if a sequence of finite graphs converges to a hyperfinite limit, then the sequence itself is hyperfinite.
\end{abstract}


\section{Introduction}\label{intro}

While studying asymptotic properties of finite graphs, it can be worthwhile to introduce
compactifications of various collections of finite graphs.
In one such compactification, which is particularly suited to bounded degree graphs,
the notion of convergence is the local weak convergence, as studied in~\cite{\BSplanarLimits}
and~\cite{\AldousLyons}.
A different notion of convergence~\cite{\LovaszSzegedy} is suitable for dense graphs,
in which the number of edges is roughly proportional to $|V(G)|^2$.
Here, we will focus on the bounded degree case.

The purpose of this note is to extend the study of the local weak convergence of graphs, whose
definition will be recalled below, by noting its properties in relation to the collection
of hyperfinite graphs.
A collection of finite graphs $\mathcal G$ is {\bf $(k,\eps)$-hyperfinite} if
every $G\in\mathcal G$ has a set of edges $S\subset E(G)$
such that $|S|\le\eps\,|V(G)|$ and every connected component of $G\setminus S$ has
at most $k$ vertices.
We say that $\mathcal G$ is {\bf hyperfinite}, if for every $\eps>0$ there is some finite $k=k(\eps)$
such that $\mathcal G$ is $(k,\eps)$-hyperfinite.
Many interesting collections of graphs are hyperfinite. 
For example, as noted in~\cite[Theorem 3]{\LiptonTarjanApplications},
it follows from the planar separator theorem that the set of planar graphs with maximal
degree at most $M$ is hyperfinite for every $M<\infty$.

Although the notion of hyperfiniteness appears implicitly in the literature
(i.e.,~\cite{\LiptonTarjanApplications}), as far as we know
G\'abor Elek~\cite{\ElekCost,\ElekLtwo,\ElekRegularity} was the first to give it a name and
propose its systematic study.

We now prepare to introduce the notion of local weak convergence.
A {\bf rooted graph} is a pair $(G,o)$, where $o\in V(G)$.
An {\bf isomorphism} of rooted graph $\phi:(G,o)\to(G',o')$ is
an isomorphism of the underlying graphs which satisfies $\phi(o)=o'$.
A graph is {\bf locally-finite} if each vertex is incident with only finitely many edges.
Let $\G$ denote the collection of all isomorphism classes of connected,
locally-finite rooted graphs, and for $M\in\N_+$ let $\GG\subset\G$
denote the subcollection consisting of rooted graphs with maximal degree at most $M$.
For $(G,o)\in\G$ and $r\ge 0$ let $B_G(o,r)$ denote
the subgraph of $G$ spanned by the vertices at distance at most $r$ from $o$.
If $(G,o),(G',o')\in\G$ and $r$ is the largest integer such that 
$\bigl(B_G(o,r),o\bigr)$ is rooted-graph  isomorphic to $\bigl(B_{G'}(o',r),o'\bigr)$,
then set $\rho\bigl((G,o),(G',o')\bigr)=1/r$, say.
Also take $\rho\bigl((G,o),(G,o)\bigr)=0$.
Then $\rho$ is metric on $\G$.
Let $\MM$ [respectively, $\M$] denote the space of all probability measures on $\GG$ [resp.\ $\G$]
that are measurable with respect to the Borel $\sigma$-field of $\rho$.
Then $\MM$ is endowed with the topology of weak convergence, and is compact in this
topology.

Let $\GfinM$ denote the collection of all isomorphism classes of finite graphs
with maximal degree at most $M$, and let $G\in\GfinM$.
Let $o\in V(G)$ be chosen
randomly-uniformly, and let $G_o$ be the connected component
of $o$. Then the law of (the isomorphism type) of $(G_o,o)$
is an element of $\MM$, which will be denoted by $\Psi(G)$.
This defines a mapping $\Psi:\GfinM\to\MM$.
The restriction of this mapping to the set of connected graphs is
injective.

Let $\MM^0$ denote the closure of $\Psi(\GfinM)$ in $\MM$.
It has been observed in~\cite{\BSplanarLimits} that every $\mu\in\MM^0$
satisfies an intrinsic version of the Mass Transport Principle (MTP).
The MTP has been invented by Olle H\"aggstr\"om in the setting of
regular trees~\cite{MR98f:60207}, and used extensively by Benjamini, Lyons, Peres and Schramm
in the more general setting of transitive and quasi-transitive graphs~\cite{MR99m:60149}.
Since then, there have been numerous successful applications of the MTP.
A measure $\mu\in\M$ satisfies the intrinsic MTP (iMTP) if
for every non-negative function $f(G,x,y)$, which takes as arguments
a connected graph and two vertices in it, and depends only on the isomorphism
type of the triple $(G,x,y)$, we have
$$
\int \sum_{x\in V(G)} f(G,x,o)\,d\mu(G,o) =
\int \sum_{y\in V(G)} f(G,o,y)\,d\mu(G,o)\,.
$$
Let $\UU$ denote the set of $\mu\in\M$ satisfying the iMTP and set $\UU_M:=\UU\cap\MM$.
The easy proof that $\MM^0\subset\UU$ proceeds by first
noting that $\Psi(\GfinM)\subset \UU$, and then observing that $\UU_M$
is closed in $\MM$.

Aldous and Lyons~\cite{\AldousLyons} coined the term {\bf unimodular} for
measures in $\UU$, and studied the properties of such measures. They also raised
the fundametal problem to determine if $\MM^0=\UU_M$; that is,
whether every $\mu\in\UU_M$ can be approximated by elements of $\Psi(\GfinM)$.
\medskip

If $\nu$ is a measure on triples $(G,o,S)$, where $G$ is a connected
graph, $o\in V(G)$, and $S$ is some structure on $G$ (such as a labeling
of the edges or vertices, or a subgraph, etc), then $\nu$ is {\bf unimodular}
if it satisfies the above iMTP, where $f$ is allowed to depend on $S$.

A measure $\mu\in\UU$ is {\bf $(k,\eps)$-hyperfinite} if there is a measure
$\nu$ on triples $(G,o,S)$, where the projection of $\nu$ to $(G,o)$
is $\mu$, $S\subset E(G)$, every connected component of $G\setminus S$ has
at most $k$ vertices, $\nu$ is unimodular and the $\nu$-expected number of
edges in $S$ adjacent to $o$ is at most $2\,\eps$.

The reason for choosing $2\,\eps$ instead of $\eps$ is the following.
If $G_0$ is a finite finite graph that is $(k,\eps)$-hyperfinite,
and $\mu=\Psi(G_0)$, then $\mu$ is $(k,\eps)$-hyperfinite.
This is because if $S_0\subset E(G_0)$ satsifies $|S_0|\le\eps |V(G_0)|$,
then the uniform law on the triples $(G_0,o,S_0)$ where $o\in V(G_0)$
is unimodular and the expected number of edges of $S_0$ adjacent to $o$
is $2\,|S_0|/|V(G)|\le2\,\eps$.
Conversely we show in Lemma~\ref{l.warm}
below that if $G\in \GfinM$ and $\Psi(G)$ is
$(k,\eps)$-hyperfinite, then $G$ is also $(k,\eps)$-hyperfinite.

A measure $\mu\in\UU$ is {\bf hyperfinite} if for every $\eps>0$ there is a finite $k=k(\eps)$
such that $\mu$ is $(k,\eps)$-hyperfinite.
\medskip

In this paper, we prove the following theorem.

\begin{theorem}\label{t.1} Let $M\in\N_+$,
and let $G_1,G_2,\dots$ be a sequence in $\GfinM$ with $\Psi(G_j)$ converging weakly
to some $\mu\in\MM$. Then $\mu$ is hyperfinite if and only if the sequence
$\{G_1, G_2,\dots\}$ is hyperfinite.
\end{theorem}

In the forthcoming paper~\cite{\BSSpropertyTesting},
we give an application of this theorem to planarity testing.

The easy direction of the theorem is to show that $\mu$ is hyperfinite if
$\{G_1,G_2,\dots\}$ is hyperfinite.
The opposite direction is an immediate consequence of the following,
more quantitive, version.

\begin{theorem}\label{t.quant}
Suppose that $\mu\in\UU_M$ is $(k,\eps)$-hyperfinite,
where $0<\eps<1$ and $k<\infty$.
Then there is some open neighborhood of $\mu$ such that
every $\mu'\in\UU_M$ within that neighborhood is
$(k,\tilde \eps)$-hyperfinite, where
$$
\tilde \eps:= 3\,\log(2\,M/\eps)\,\eps\,.
$$
\end{theorem}

Theorem~\ref{t.quant} is quantitative in the sense that the dependence
of $\tilde\eps$ on $\eps$ is explicit. 
In~\cite{\BSSpropertyTesting},
we present a finitary variant of Theorem~\ref{t.quant},
which bounds the size of the neighborhood. The proof there is
also a finitary variant of the proof below.

I.~Benjamini (private communication) points out that hyperfiniteness of
measures in $\UU_M$ is an appropriate analogue for amenability in the
setting of transitive graphs. In particular, the Burton-Keane~\cite{MR90g:60090}
argument can show that percolation on a sample from $\UU_M$ can
produce at most one infinite cluster a.s.\ 
Extending a conjecture from~\cite{\BSpyond} in the transitive setting,
Benjamini also conjectures
the converse: that if $\mu\in\UU_M$ is not hyperfinite, then there is
some parameter $p\in[0,1]$ such that Bernoulli($p$) percolation
on the sample from $\mu$ has more than one infinite cluster with positive probability.
See~\cite{2003math......6355A} for a related discussion in the context of percolation
on finite graphs and in particular the discrete hypercube.

\section{Proofs}

As a warm up, we prove the following.

\begin{lemma}\label{l.warm}
Let $\eps>0$, $k\in\N_+$ and $G\in\GfinM$.
Then $G$ is $(k,\eps)$-hyperfinite
if and only if $\Psi(G)$ is $(k,\eps)$-hyperfinite.
\end{lemma}
\proof
The \lq\lq only if\rq\rq\ direction was sketched in the introduction,
and hence we presently only prove the \lq\lq if\rq\rq\ direction.
Assume that $\Psi(G)$ is $(k,\eps)$-hyperfinite.
Then there is a unimodular probability measure $\nu$
on triples $(G',S,o)$ such that $G'$ is $\nu$-a.s.\ 
a connected component
of $G$, $S\subset E(G')$, $o\in V(G')$, the $\nu$-distribution
of $o$ is uniform in $V(G)$, the connected components of $G'\setminus S$
all have at most $k$ vertices $\nu$-a.s.\ and the $\nu$-expected number of
edges in $S$ that are adjacent to $o$ is at most $2\,\eps$.

Let $S(v)$ denote the edges of $S$ incident to a vertex $v$.
Let $G_1,G_2,\dots,G_m$ be the connected components of $G$.
Let $\ev A_i$ be the event that $G'$ is isomorphic to $G_i$.
Define $g(G',S,x,y)$ as the number of edges of $S$ that
are incident to $x$ and $f_i(G',S,x,y)=1_{\ev A_i}\,g(G',S,x,y)$.
Since $\nu$ is unimodular, the iMTP may be applied to $f_i$:
\begin{equation}
\label{e.im1}
\int\sum_{y\in V(G')} f_i(G',S,o,y)\,d\nu
=\int\sum_{x\in V(G')} f_i(G',S,x,o)\,d\nu
\,.
\end{equation}
The left hand side of~\eqref{e.im1} is equal to
$$
\nu\bigl[|S(o)|\,|V(G')|\,1_{\ev A_i}\bigr]=
 |V(G_i)|\,\nu\bigl[|S(o)|\,1_{\ev A_i}\bigr],
$$
while the right hand side is
$ \nu\bigl[2\,|S|\,1_{\ev A_i}\bigr]$. 
Hence,
$$
 |V(G_i)|\,\nu\bigl[|S(o)|\,1_{\ev A_i}\bigr]=
 \nu\bigl[2\,|S|\,1_{\ev A_i}\bigr].
$$
Therefore, there is some $S_i\subset E(G_i)$ such
that each connected component of $G_i\setminus S_i$
has at most $k$ vertices and
\begin{equation}
\label{e.Si}
|S_i|\,\nu[\ev A_i]
\le \frac{|V(G_i)|}{2}\,\nu\bigl[ |S(o)|\,1_{\ev A_i}\bigr].
\end{equation}
Define $\tilde S:=\bigcup_{i=1}^m S_i$.
Note that $\nu[\ev A_i]\,|V(G)|$ is the number of vertices
of $G$ that are in connected components of $G$ that are
isomorphic to $G_i$ and hence
$t_i:=\nu[\ev A_i]\,|V(G)|/|V(G_i)|$
is the number of connected
components of $G$ that are isomorphic to $G_i$.
Note that
$$
\sum_{i=1}^m \frac{1_{\ev A_i}}{t_i}=1\,.
$$
 Dividing~\eqref{e.Si} by
$\nu[\ev A_i]$ and summing over $i$ gives
$$
|\tilde S| \le \frac{|V(G)|}2\,\nu\Bigl[|S(o)|\sum_i \frac {1_{\ev A_i}}{t_i}
\Bigr]=
\frac{|V(G)|}2\,\nu\bigl[|S(o)|\bigr]\le 2\,\eps\,|V(G)|\,,
$$
which completes the proof.
\QED

The idea of the proof of Theorem~\ref{t.quant} is to replace
the random set of edges $S$ by a set of edges that still separates $G$ into connected
components of size at most $k$, but has the extra feature that for any $e\in E(G)$
the event $e\in S$ depends only on the local structure of $G$ near $e$ and some randomness.
Thus, it is easy to adapt the law of this new $S$ to every random $G'$ sufficiently close to $G$.

The proof of the theorem is complicated by a completely uninteresting point 
of a technical nature.
For the sake of simplicity, we choose to 
present a proof that is not entirely precise, but does convey
the essential ideas, and can be adapted to be completely correct.
(See,~\cite{\BSS}.) 

\proofof{Theorem~\ref{t.quant}}
Fix some $\eps_0\in(0,1/2)$.
Let $\P$ be a unimodular measure on triples $(G,S,o)$ such that the marginal
law of $(G,o)$ is $\mu$ the $\P$-expected number of edges
of $S$ incident with $o$ is at most $2\,\eps$, and every connected component of
$G\setminus S$ has at most $k$ vertices.

For any $K\subset V(G)$ and any $v\in K$, let $p(K):= \Pb{K(v)=K\md (G,o)}$.
(Here, we are not entirely precise, since the expression $\Pb{K(v)=K\md (G,o)}$
does not have the standard meaning of the probability of an event conditional
on a $\sigma$-field. This is because $K$ and $v$ do not have an a priori meaning without $G$. This
difficulty is not too hard to overcome, at the expense of obfuscating the proof.)
Clearly, $p(K)$ does not depend on the choice of the particular $v\in K$.
Also, a simple argument using the iMTP shows that $p(K)$ is not effected by a change of the
basepoint $o$; in other words, $p(K)$ is really a function of the isomorphism
type of the pair $(K,G)$.

Given $(G,o)$, and $v\in V(G)$, let $\mathfrak K_v$
denote the collection of all connected sets of vertices
of cardinality at most $k$ that contain $v$.
Obviously, $|\mathfrak K_v|$ is bounded by a function of $k$ and $M$.
For all $r\in\N_+$ let $\mB_r$
denote the (isomorphism type of the) rooted graph $\bigl(B(o,r),o\bigr)$.
For $r>k$, and $K\in\mathfrak K_o$, we set
$$
p_r(K):= \Pb{K(o)=K\md \mB_r}.
$$
Clearly, $p_r(K)$ is a bounded martingale with respect to $r$
and $\lim_{r\to\infty} p_r(K)=p(K)$ a.s.
Consequently, there is some finite $R$ such that
$$
\E \sum_{K\in\mathfrak K_o} |p_R(K)-p(K)|<\eps_0\,.
$$
Fix such an $R$.

For $K\in \mathfrak K_o$ let
$$
N_R(K):=\bigcup_{v\in K} B(v,R)\,,
\qquad
\tilde p(K)=\tilde p(K,G):= \Pb{K=K(o)\md N_R(K)}.
$$
By unimodularity,  $\tilde p(K)$ does not depend on the choice of basepoint $o$ in $K$.
Therefore, by changing basepoint this also defines $\tilde p(K)$ for all
$\mathfrak K_v$, $v\in V(G)$.
Note that $\tilde p$ is a function whose arguments
are a graph and a connected set of size at most $k$ in it,
and $\tilde p$ depends only on the $R$-neighborhood of the set.
However, $\tilde p$ is only defined if its arguments can occur
in $(G,o)$. We extend the definition of $\tilde p$ to such arguments that
have probability $0$ of occuring in $(G,o)$ (if such exist)
by setting it equal to $1$ on such arguments.

Let $\mathfrak K:=\bigcup_{v\in V(G)} \mathfrak K_v$.
Given $(G,o)$, let $(X_K:K\in\mathfrak K)$ be independent random variables
such that 
$$
\Pb{X_K=1\md (G,o)}=\min\{2\,\log (1/\eps_0)\,\tilde p(K),1\}\,,
$$
and $X_K\in\{0,1\}$ a.s.
Let $F:=\bigcup_{K:X_K=1} \p K$,
where $\p K$ denotes the edge-boundary of $K$,
and let $W:=\bigcup\{K\in\mathfrak K:X_K=1\}$.
Let $\tilde F$ denote the set of all edges incident with some vertex in $V(G)\setminus W$,
and set $S':=F\cup\tilde F$. It is easy to verify that the law of
$(G,S',o)$ is unimodular.
Also,
it is obvious that every connected component of $G\setminus S'$ has cardinality at most $k$.

\smallskip

Our present goal is to estimate $\Eb{|E(o)\cap S'|}$. We start by
estimating $\Eb{|E(o)\cap \tilde F|}$.
Let $\ev A$ denote the event $\sum_{K\in \mathfrak K_o} \tilde p(K) < 1/2$.

\begin{lemma}\label{sump}
$$
\Pb{\ev A} < 2\,\eps_0\,.
$$
\end{lemma}
\proof
Fix some $K\in\mathfrak K_o$.
Since the sequence $p_R(K),\tilde p(K), p(K)$ is a martingale given $B(o,k)$,
we have
$$
\Eb{|\tilde p(K)-p(K)|\md B(o,k)}\le \Eb{|p_R(K) -p(K)|\md B(o,k)}\,.
$$
Our choice of $R$ therefore gives
$$
\EB{ \sum_{K\in\mathfrak K_o} |\tilde p(K)-p(K)|} <\eps_0\,.
$$
Since $\sum_{K\in\mathfrak K_o} p(K)=1$, we get
$$
\EB{\Bl|1- \sum_{K\in\mathfrak K_o} \tilde p(K)\Br|} <\eps_0\,,
$$
and Markov's inequality completes the proof of the lemma.
\QED

\begin{lemma}
$$
\Pb{o\notin W} < 3\,\eps_0\,.
$$
\end{lemma}
\proof
Given $(G,o)$, the random variable $X_K$ are independent and satisfy
$\Pb{X_K=1\md (G,o)}=\min\bl\{1,2\,\log(1/\eps_0)\,\tilde p(K)\br\}$.
Hence,
\begin{align*}
\Pb{o\notin W\md (G,o)}
&
=
\prod_{K\in\mathfrak K_o}
\Bl( 1- \Pb{X_K=1\md (G,o)}\Br)
\\&
\le
\exp\Bl( -\sum_{K\in\mathfrak K_o} 2\,\log(1/\eps_0)\,\tilde p(K)\Br)
\\&
\le \exp\bl(-2\,\log(1/\eps_0)\,(1/2)\br)+1_{\ev A}
=\eps_0+1_{\ev A}\,.
\end{align*}
Now Lemma~\ref{sump} completes the proof.
\QED

The lemma implies
\begin{equation*}
\Eb{|E(o)\cap \tilde F|}\le 3\,M\,\eps_0\,,
\end{equation*}
where $M$ is the bound on the degrees in $G$.
  
Now,
\begin{align*}
\Eb{|E(o)\cap F|\md (G,o)}
&
\le 
\sum_{K\in\mathfrak K_o}
|\p K\cap E(o)|\,
\Pb{X_K=1\md (G,o)}
\\&
\le
2\,
\log(1/\eps_0)
\sum_{K\in\mathfrak K_o}
|\p K\cap E(o)|\,
\tilde p(K)\,.
\end{align*}
Since,
$$
\EB{
\sum_{K\in\mathfrak K_o}
|\p K\cap E(o)|\,
\tilde p(K)}=\Eb{|E(o)\cap \p K(o)|}=
\Eb{E(o)\cap S}\le2\,\eps\,,
$$
we get
$$
\Eb{|E(o)\cap F|}\le 4\,\log(1/\eps_0)\,\eps\,.
$$
Thus, we have,
$$
\Eb{|E(o)\cap S'|}\le 4\,\log(1/\eps_0)\,\eps+ 3\,\eps_0\,M\,.
$$
We now choose $\eps_0:=\eps/(2\,M)$.

Observe that the set $S'$ is chosen in a very local way.
Namely, given $(G,o)$ you can decide if $e\in S'$ based purely on a fixed radius
neighborhood of $e$ and some coin flips that are associated with this neighborhood.
Given another unimodular $(\tilde G,\tilde o)$, we 
can choose a $\tilde S\subset E(\tilde G)$ according to the same procedure as $S'$ is obtained
for $G$.
If $(\tilde G,\tilde o)$ is sufficiently close to $(G,o)$, then the expected size
of the set of edges in $\tilde S$ adjacent to $\tilde o$ will be close to the
corresponding quantity in $G$.
This proves the theorem.
\QED

\proofof{Theorem~\ref{t.1}}
It is easy to verify that the set of $\nu\in\UU_M$ that are $(k,\eps)$-hyperfinite
is closed. Therefore, Lemma~\ref{l.warm} implies that $\mu$ is hyperfinite
if $\{G_1, G_2,\dots\}$ is hyperfinite. The converse follows immediately from
Theorem~\ref{t.quant} and Lemma~\ref{l.warm}.
\QED

\bigskip
\noindent
{\bf Acknowledgments}.
We are grateful to G\'abor Elek, Christian Borgs and Russell Lyons for beneficial conversations.

\bibliographystyle{halpha}
\bibliography{mr,prep,notmr}

\def\polhk#1{\setbox0=\hbox{#1}{\ooalign{\hidewidth
  \lower1.5ex\hbox{`}\hidewidth\crcr\unhbox0}}} \def\cprime{$'$}
  \def\cprime{$'$}
\begin{thebibliography}{{Ele}07b}

\bibitem[AB03]{2003math......6355A}
O.~{Angel} and I.~{Benjamini}.
\newblock {A Phase Transition for the Metric Distortion of Percolation on the
  Hypercube}, 2003, arXiv:math/0306355.

\bibitem[AL06]{2006math......3062A}
D.~{Aldous} and R.~{Lyons}.
\newblock {Processes on Unimodular Random Networks}, 2006, arXiv:math/0603062.

\bibitem[BK89]{MR90g:60090}
R.~M. Burton and M.~Keane.
\newblock Density and uniqueness in percolation.
\newblock {\em Comm. Math. Phys.}, 121(3):501--505, 1989.

\bibitem[BLPS99]{MR99m:60149}
I.~Benjamini, R.~Lyons, Y.~Peres, and O.~Schramm.
\newblock Group-invariant percolation on graphs.
\newblock {\em Geom. Funct. Anal.}, 9(1):29--66, 1999.

\bibitem[BS96]{MR97j:60179}
Itai Benjamini and Oded Schramm.
\newblock Percolation beyond $\bold {Z}\sp d$, many questions and a few
  answers.
\newblock {\em Electron. Comm. Probab.}, 1:no.\ 8, 71--82 (electronic), 1996.

\bibitem[BS01]{MR1873300}
Itai Benjamini and Oded Schramm.
\newblock Recurrence of distributional limits of finite planar graphs.
\newblock {\em Electron. J. Probab.}, 6:no. 23, 13 pp. (electronic), 2001.

\bibitem[BSS07]{BSSpropertyTesting}
I.~Benjamini, O.~Schramm, and A.~Shapira.
\newblock Every minor-closed property of sparse graphs is testable, 2007.
\newblock Preprint.

\bibitem[{Ele}06]{2006math......8474E}
G.~{Elek}.
\newblock {The combinatorial cost}, 2006, arXiv:math/0608474.

\bibitem[{Ele}07a]{2007arXiv0709.1261E}
G.~{Elek}.
\newblock {$L^2$-spectral invariants and convergent sequences of finite
  graphs}, 2007, arXiv:0709.1261.

\bibitem[{Ele}07b]{arXiv:0711.2800}
G.~{Elek}.
\newblock A regularity lemma for bounded degree graphs and its applications:
  Parameter testing and infinite volume limits, 2007, arXiv:0711.2800.

\bibitem[H{\"a}g97]{MR98f:60207}
Olle H{\"a}ggstr{\"o}m.
\newblock Infinite clusters in dependent automorphism invariant percolation on
  trees.
\newblock {\em Ann. Probab.}, 25(3):1423--1436, 1997.

\bibitem[LS06]{MR2274085}
L{\'a}szl{\'o} Lov{\'a}sz and Bal{\'a}zs Szegedy.
\newblock Limits of dense graph sequences.
\newblock {\em J. Combin. Theory Ser. B}, 96(6):933--957, 2006.

\bibitem[LT80]{MR584516}
Richard~J. Lipton and Robert~Endre Tarjan.
\newblock Applications of a planar separator theorem.
\newblock {\em SIAM J. Comput.}, 9(3):615--627, 1980.

\end{thebibliography}

\end{document}